December 18, 2017

# POWER SERIES WITH INVERSE BINOMIAL COEFFICIENTS AND HARMONIC NUMBERS

KHRISTO N. BOYADZHIEV
Ohio Northern University
Department of Mathematics and Statistics
Ada, Ohio 45810, USA
k-boyadzhiev@onu.edu

**ABSTRACT** In this note we compute the generating function for the numbers $\binom{2n}{n}^{-1} H_n$ in terms of elementary functions and dilogarithms.

## 1. Introduction and main result

Let $\binom{2n}{n} = \frac{(2n)!}{(n!)^2}$ be the central binomial coefficients. Series with these numbers and their reciprocals $\binom{2n}{n}^{-1}$ have gained considerable attention in mathematics. Many mathematicians became interested in series involving $\binom{2n}{n}^{-1}$ after Roger Apéry [1] used the representation

$$\zeta(3) = \frac{5}{2} \sum_{n=1}^{\infty} \frac{(-1)^{n-1}}{n^3} \binom{2n}{n}^{-1}$$

to show that $\zeta(3)$ is an irrational number. Interesting contributions were made later by Lehmer [9] and also Glasser [6], who worked with the power series $\sum_{n=1}^{\infty}\binom{2n}{n}^{-1} n^m z^n$. The attention toward this topic was reinforced when it was realized that power series of the form

$$\sum_{n=1}^{\infty}\binom{2n}{n}^{-1}\frac{z^n}{n^k} \text{ and also } \sum_{n=1}^{\infty}\binom{2n}{n}^{-1} S_p \frac{z^n}{n}$$

(here $S_p$ are certain harmonic sums) were closely related to the computation of Feynman diagrams (see [4], [7], [8]). In these papers various series with the above structure were computed in closed form in terms of certain special functions.

As usual, let $H_n = 1 + \frac{1}{2} + \ldots + \frac{1}{n}$, $H_0 = 0$, be the harmonic numbers. Using a special series transformation formula the author computed in [3] the generating function for the numbers $\binom{2n}{n} H_n$ in terms of elementary functions. Using another series transformation formula (Lemma 2 below) we compute here the generating function for the numbers $\binom{2n}{n}^{-1} H_n$ explicitly in terms of elementary functions and dilogarithms. The form of this generating function makes it possible to consider analytic continuation beyond the disk of convergence.

Our main result is the theorem:

**Theorem 1.** *For $|z|<1$ we have:*

*Part 1. With $t = \dfrac{z}{1+z}, |t| \leq 1$*

$$\sum_{n=0}^{\infty}\binom{2n}{n}^{-1} H_n (-1)^n 4^n z^n \qquad (1)$$

$$= (1-t)\left[\sqrt{t}\log\frac{1-\sqrt{t}}{1+\sqrt{t}}\left(1+\frac{1}{4}\log\frac{1-t}{4}\right) - \frac{\sqrt{t}}{2}\left(\text{Li}_2\left(\frac{1+\sqrt{t}}{2}\right) - \text{Li}_2\left(\frac{1-\sqrt{t}}{2}\right)\right)\right]$$

$$= -2z + 4z^2 - \frac{88}{15}z^3 + \frac{160}{21}z^4 - \frac{8768}{945}z^5 + \ldots \ .$$

*Part 2. With $t = \dfrac{z}{1-z}, |t| \leq 1$*

$$\sum_{n=0}^{\infty} \binom{2n}{n}^{-1} H_n \, 4^n z^n \qquad (2)$$

$$= (1+t)\left[\left(2+\frac{1}{2}\log\frac{1+t}{4}\right)\sqrt{t}\arctan\sqrt{t} - \frac{\sqrt{-t}}{2}\left(\operatorname{Li}_2\left(\frac{1+\sqrt{-t}}{2}\right) - \operatorname{Li}_2\left(\frac{1-\sqrt{-t}}{2}\right)\right)\right]$$

$$= 2z + 4z^2 + \frac{88}{15}z^3 + \frac{160}{21}z^4 + \frac{8768}{945}z^5 + \ldots \quad .$$

Here, of course,

$$\operatorname{Li}_2(x) = \sum_{n=1}^{\infty} \frac{x^n}{n^2}$$

is the dilogarithm. The proof of the theorem is in Section 2. Concerning the functions with $\sqrt{t}$ and $\sqrt{-t}$ see the note after Lemma 3 below.

If we denote by $F(t)$ the function in the brackets on the RHS in (1), then the generating function in (1) can be written shortly

$$G(z) = \sum_{n=0}^{\infty} \binom{2n}{n}^{-1} H_n (-1)^n 4^n z^n = \frac{1}{1+z} F\left(\frac{z}{1+z}\right).$$

Correspondingly, (2) becomes

$$G(-z) = \sum_{n=0}^{\infty} \binom{2n}{n}^{-1} H_n 4^n z^n = \frac{1}{1-z} F\left(\frac{-z}{1-z}\right).$$

Note that the radius of convergence of the above series is $R=1$.

**Remark.** Setting $z = \frac{1}{4}$ in (2) we compute $t = \frac{1}{3}$ and

$$\sum_{n=0}^{\infty} \binom{2n}{n}^{-1} H_n = (4-\log 3)\frac{\pi}{9\sqrt{3}} - \frac{2i}{3\sqrt{3}}\left(\operatorname{Li}_2\left(\frac{1}{2}+\frac{i}{2\sqrt{3}}\right) - \operatorname{Li}_2\left(\frac{1}{2}-\frac{i}{2\sqrt{3}}\right)\right)$$

which is exactly the evaluation obtained by M. Genchev [5, p.219, (9)]. With $z = \frac{1}{4}$ and $t = \frac{1}{5}$ in (1) we also compute the alternating series

$$\sum_{n=0}^{\infty}\binom{2n}{n}^{-1}(-1)^n H_n = \frac{4-\log 5}{5\sqrt{5}}\log\frac{\sqrt{5}-1}{\sqrt{5}+1} - \frac{1}{2\sqrt{5}}\left(\operatorname{Li}_2\left(\frac{\sqrt{5}+1}{2\sqrt{5}}\right) - \operatorname{Li}_2\left(\frac{\sqrt{5}-1}{2\sqrt{5}}\right)\right).$$

The proof of the theorem is given in the next section. It is based on two lemmas, the first of which is an Euler type series transformation formula.

**Lemma 2.** *Let*

$$f(z) = \sum_{n=0}^{\infty} a_n z^n$$

*be a function analytic in a neighborhood of $z=0$. Then for any $z$ with $|z|$ small enough to assure convergence we have*

$$\sum_{n=0}^{\infty} a_n H_n z^n + \log(1+z) f(z) = \frac{1}{1+z}\sum_{n=0}^{\infty}\left(\frac{z}{z+1}\right)^n H_n \left\{\sum_{k=0}^{n}\binom{n}{k} a_k\right\}. \tag{3}$$

The proof of this lemma can be found in [2].

The next lemma gives the generating function of the numbers $\dfrac{H_n}{2n-1}$ in explicit form.

**Lemma 3.** *For every $|t|<1$ we have*

*Part 1.*

$$\sum_{n=1}^{\infty} H_n \frac{t^n}{2n-1} = \log(1-t) + \left(-1+\log\sqrt{2}+\frac{1}{4}\log(1-t)\right)\sqrt{t}\,\log\frac{1-\sqrt{t}}{1+\sqrt{t}}$$

$$+ \frac{\sqrt{t}}{2}\left[\operatorname{Li}_2\left(\frac{1+\sqrt{t}}{2}\right) - \operatorname{Li}_2\left(\frac{1-\sqrt{t}}{2}\right)\right] \tag{4}$$

*Part 2.*

$$\sum_{n=1}^{\infty} H_n \frac{(-1)^n t^n}{2n-1} = \log(1+t) + \left(-2+\ln 2+\frac{1}{2}\log(1+t)\right)\sqrt{t}\,\arctan\sqrt{t}$$

$$+ \frac{\sqrt{-t}}{2}\left[\operatorname{Li}_2\left(\frac{1+\sqrt{-t}}{2}\right) - \operatorname{Li}_2\left(\frac{1-\sqrt{-t}}{2}\right)\right]. \tag{5}$$

Obviously, the series in (4) converges also for $t = -1$, while the series in (5) converges also for $t = 1$. In view of its technical character, the proof is given in the Appendix. The expressions with $\sqrt{t}$ or $\sqrt{-t}$ do not cause trouble. For instance,

$$\sqrt{t}\arctan\sqrt{t} = \sum_{n=1}^{\infty}\frac{(-1)^{n-1}t^n}{2n-1}$$

is a well-defined analytic function in the unit disk and replacing $t$ by $-t$ is legitimate, as we can write $\sqrt{-t} = i\sqrt{t}$ when $t \geq 0$. See also equations (7) and (11) below. The same is true for the last terms in (1) and (2), as the function is analytic in the unit disk

$$\frac{\sqrt{t}}{2}\left[\operatorname{Li}_2\left(\frac{1+\sqrt{t}}{2}\right) - \operatorname{Li}_2\left(\frac{1-\sqrt{t}}{2}\right)\right] = (\ln 4)t + \left(\frac{-1}{3} + \frac{\ln 4}{3}\right)t^2 + \left(\frac{-7}{30} + \frac{\ln 4}{5}\right)t^3\ldots\ .$$

## 2. Proof of the theorem

First we prove Part 1. In Lemma 2 we set

$$a_n = \binom{2n}{n}^{-1}(-1)^n 4^n$$

The generating function of these numbers is known and can be written in the form

$$f(z) = \sum_{n=0}^{\infty}\binom{2n}{n}^{-1}(-1)^n 4^n z^n = \frac{1}{1+z}\left[1 + \frac{1}{2}\sqrt{\frac{z}{1+z}}\log\frac{1-\sqrt{\frac{z}{1+z}}}{1+\sqrt{\frac{z}{1+z}}}\right] \qquad (6)$$

This is a modification of Sprugnoli's result [10, Theorem 2.1] (see (10) below) obtained by using the identity

$$\sqrt{-t}\arctan\sqrt{-t} = \frac{\sqrt{t}}{2}\log\frac{1-\sqrt{t}}{1+\sqrt{t}}\ . \qquad (7)$$

It is also known that ([10, Theorem 4.5])

$$\sum_{k=0}^{n}\binom{n}{k}a_k = \sum_{k=0}^{n}\binom{n}{k}\binom{2n}{n}^{-1}(-1)^k 4^k = \frac{1}{1-2n} \qquad (8)$$

Therefore, equation (3) takes the form

$$\sum_{n=0}^{\infty} \binom{2n}{n}^{-1} H_n (-1)^n 4^n z^n + \log(1+z) f(z) = \frac{1}{1+z} \sum_{n=0}^{\infty} \left(\frac{z}{z+1}\right)^n \frac{H_n}{1-2n}. \qquad (9)$$

With the substitution

$$t = \frac{z}{1+z}, \quad z = \frac{t}{1-t}, \quad \frac{1}{1+z} = 1-t$$

we obtain from (6)

$$f(z) = (1-t)\left[1 + \frac{1}{2}\sqrt{t} \log \frac{1-\sqrt{t}}{1+\sqrt{t}}\right], \quad \log(1+z) = -\log(1-t)$$

and (9) can be written as

$$\sum_{n=0}^{\infty} \binom{2n}{n}^{-1} H_n (-1)^n 4^n z^n = (1-t)\left[\log(1-t)\left(1 + \frac{1}{2}\sqrt{t} \log \frac{1-\sqrt{t}}{1+\sqrt{t}}\right) + \sum_{n=0}^{\infty} \frac{H_n}{1-2n} t^n\right].$$

All we need to do now in order to finish the proof is to use Part 1 of Lemma 3 (i.e. equation (4)) for the sum on the right hand side and then simplify the resulting expression.

Part 2 of the theorem is obtained from Part 1 by changing $z$ to $-z$ and using Part 2 of Lemma 3. We also need some small modifications. Namely, instead of using the representation (6) we use the representation (see [10, Theorem 2.1])

$$f(-z) = \sum_{n=0}^{\infty} \binom{2n}{n}^{-1} 4^n z^n = \frac{1}{1-z}\left(1 + \sqrt{\frac{z}{1-z}} \arctan \sqrt{\frac{z}{1-z}}\right) \qquad (10)$$

and the substitution

$$t = \frac{z}{1-z}, \quad z = \frac{t}{1+t}, \quad \frac{1}{1-z} = 1+t \ .$$

Also, we use the identity

$$\sqrt{-t} \log \frac{1-\sqrt{-t}}{1+\sqrt{-t}} = 2\sqrt{t} \arctan \sqrt{t} \ , \qquad (11)$$

connecting (6) and (10). Details are left to the reader. The theorem is proved.

## 3. Appendix

Here we prove Lemma 3, Part 1. The second part follows from the first by the substitution $t \to -t$ and by using identity (11). Let now $|t|<1$ and set

$$y(t) = \sum_{n=1}^{\infty} H_n \frac{t^n}{1-2n} \quad .$$

Then

$$\frac{d}{dt}\left(\frac{y(t^2)}{t}\right) = -\sum_{n=1}^{\infty} H_n t^{2n-2} = \frac{\log(1-t^2)}{t^2(1-t^2)} = \left(\frac{1}{t^2} + \frac{1}{1-t^2}\right)\log(1-t^2). \quad (12)$$

We have used here the well-known generating function: For $|x|<1$,

$$-\sum_{n=1}^{\infty} H_n x^n = \frac{\log(1-x)}{1-x} \quad .$$

From (12),

$$\frac{y(t^2)}{t} = \int_0^t \frac{\log(1-x^2)}{x^2} dx + \int_0^t \frac{\log(1-x^2)}{1-x^2} dx = A + B \quad . \quad (13)$$

Integrating the first integral by parts and the second integral by using partial fractions we arrive at the indicated result (4). Here are the details:

$$A = -\int_0^t \log(1-x^2)\, d\frac{1}{x} = -\frac{\log(1-t^2)}{t} + \int_0^t \frac{2}{1-x^2} dx = -\frac{\log(1-t^2)}{t} + \log\frac{1-t}{1+t} \quad ,$$

$$B = \int_0^t \frac{\log(1-x^2)}{(1-x)(1+x)} dx = \frac{1}{2}\int_0^t \left(\frac{1}{1-x} + \frac{1}{1+x}\right)\bigl(\log(1-x) + \log(1+x)\bigr) dx$$

$$= \frac{1}{2}\left(\int_0^t \frac{\log(1-x)}{1-x} dx + \int_0^t \frac{\log(1+x)}{1+x} dx + \int_0^t \frac{\log(1+x)}{1-x} dx + \int_0^t \frac{\log(1-x)}{1+x} dx\right)$$

$$= \frac{1}{4}\log^2(1+t) - \frac{1}{4}\log^2(1-t) + \frac{1}{2}\int_{-t}^{t} \frac{\log(1+x)}{1-x} dx \quad .$$

The last integral we solve with the substitution $1-x=u$ in the following way:

$$\frac{1}{2}\int_{-t}^{t} \frac{\log(1+x)}{1-x} dx = \frac{1}{2}\int_{1-t}^{1+t} \frac{\log(2-u)}{u} du = \frac{1}{2}\int_{1-t}^{1+t} \left(\log 2 + \log\left(1 - \frac{u}{2}\right)\right)\frac{du}{u}$$

$$= \frac{1}{2}\log 2 \log\frac{1+t}{1-t} - \frac{1}{2}\sum_{n=1}^{\infty}\frac{1}{n2^n}\int_{1-t}^{1+t} u^{n-1}du = \log\sqrt{2}\log\frac{1+t}{1-t} - \frac{1}{2}\sum_{n=1}^{\infty}\frac{(1+t)^n - (1-t)^n}{n^2 2^n}.$$

The sum on the RHS is, in fact, the difference of two dilogarithms. Thus

$$B = \frac{1}{4}\log^2(1+t) - \frac{1}{4}\log^2(1-t) + \log\sqrt{2}\log\frac{1+t}{1-t} - \frac{1}{2}\left[\text{Li}_2\left(\frac{1+t}{2}\right) - \text{Li}_2\left(\frac{1-t}{2}\right)\right]$$

We can write now (using the identity $a^2 - b^2 = (a-b)(a+b)$)

$$\frac{1}{4}\log^2(1+t) - \frac{1}{4}\log^2(1-t) = -\frac{1}{4}\log(1-t^2)\log\frac{1-t}{1+t},$$

and the representation (4) follows from (13). The proof of the lemma is complete.